\documentclass[psamsfonts,english]{amsart}

    \usepackage[T1]{fontenc}
    \usepackage{babel}
    \usepackage{natbib}
    \usepackage{amsmath}
    \usepackage{amssymb}
    \usepackage{amsthm}
    \usepackage[mathcal]{eucal}
    \usepackage{url}
    \usepackage[all]{xypic}

    \newtheorem{teo}{Theorem}
    \newtheorem{cor}{Corollary}
    \theoremstyle{definition}
    \newtheorem{defn}[teo]{Definition}
    \theoremstyle{remark}
    \newtheorem{rem}[teo]{Remark}

    \newcommand{\FF}{\mathbb{F}}
    \newcommand{\QQ}{\mathbb{Q}}
    \newcommand{\ZZ}{\mathbb{Z}}
    \newcommand{\RR}{\mathbb{R}}
    
    \newcommand{\PP}{\mathbb{P}}
    \newcommand{\jac}[1]{\mathcal{J}_{#1}}
    \newcommand{\heltal}[1]{\mathfrak{O}_{#1}}

    \DeclareMathOperator{\End}{End}
    \DeclareMathOperator{\Div}{Div}
    \DeclareMathOperator{\Mat}{Mat}
    \DeclareMathOperator{\gal}{Gal}

\begin{document}

\title[Embedding Degree of Hyperelliptic Curves with CM]
{Embedding Degree of Hyperelliptic Curves with Complex Multiplication}

\author[C.R. Ravnshøj]{Christian Robenhagen Ravnshøj}

\address{Department of Mathematical Sciences \\
University of Aarhus \\
Ny Munkegade \\
Building 1530 \\
DK-8000 Aarhus C}

\email{cr@imf.au.dk}

\thanks{Research supported in part by a PhD grant from CRYPTOMAThIC}



% 14H40 Curves, Jacobians
% 14Q05 Computational aspects in algebraic geometry, curves
% 11G15 Arithmetic algebraic geometry (Diophantine geometry), Complex multiplication
% 94A60 Communication, information, Cryptography

\begin{abstract}
Consider the Jacobian of a genus two curve defined over a finite field and with complex multiplication. In this paper we show that if the \mbox{$\ell$-Sylow} subgroup of the Jacobian is not cyclic, then the embedding degree of the Jacobian with respect to $\ell$ is one.
\end{abstract}

\maketitle

\section{Introduction}

In elliptic curve cryptography it is essential to know the number of points on the curve. Cryptographically we are
interested in elliptic curves with large cyclic subgroups. Such elliptic curves can be constructed. The construction is based on the
theory of complex multiplication, studied in detail by \cite{atkin-morain}. It is referred to as the \emph{CM method}.

\cite{koblitz89} suggested the use of hyperelliptic curves to provide larger group orders. Therefore constructions of
hyperelliptic curves are interesting. The CM method for elliptic curves has been generalized to hyperelliptic curves of genus~two by
\cite{spallek}, and efficient algorithms have been proposed by \cite{weng03} and \cite{gaudry}.

Both algorithms take as input a primitive, quartic CM field $K$ (see section~\ref{sec:CMfields} for the definition of a CM field), and give as output a hyperelliptic genus~two curve $C$
defined over a prime field $\FF_p$. A prime number $p$ is chosen such that $p=x\overline x$ for a number
$x\in\heltal{K}$, where $\heltal{K}$ is the ring of integers of $K$. We have $K=\QQ(\eta)$ and
$K\cap\RR=\QQ(\sqrt{D})$, where $\eta=i\sqrt{a+b\xi}$ and
    \begin{equation*}
    \xi=\begin{cases}
    \frac{1+\sqrt{D}}{2}, & \textrm{if $D\equiv 1\mod{4}$,} \\
    \sqrt{D}, & \textrm{if $D\equiv 2,3\mod{4}$}.
    \end{cases}
    \end{equation*}
In this paper, the following theorem is established.

\begin{teo}
Let $C$ be a hyperelliptic curve of genus two defined over $\FF_p$ with $\End(C)\simeq\heltal{K}$, where $K$ is a
primitive, quartic CM field as defined in definition~\ref{def:CMfieldPrimitive}. Assume that the $p$-power Frobenius under this isomorphism is given by the number
$\omega=c_1+c_2\xi+(c_3+c_4\xi)\eta$, where $\xi$ and $\eta$ are given as above and $c_i\in\ZZ$. Consider a prime
number $\ell\mid|\jac{C}(\FF_p)|$ with $\ell\neq p$, $\ell\nmid D$ and $\ell\nmid c_2$. Assume that the $\ell$-Sylow
subgroup of $\jac{C}(\FF_p)$ is not cyclic. Then $p\equiv 1\mod{\ell}$, i.e. the embedding degree of
$\jac{C}(\FF_p)$ with respect to $\ell$ is one.
\end{teo}

\section{Hyperelliptic curves}

A hyperelliptic curve is a smooth, projective curve $C\subseteq\PP^n$ of genus at least two with a separable, degree
two morphism $\phi:C\to\PP^1$. Let $C$ be a hyperelliptic curve of genus two defined over a prime field $\FF_p$ of
characteristic~$p>2$. By the Riemann-Roch theorem there exists an embedding $\psi:C\to\PP^2$, mapping $C$ to a curve
given by an equation of the form
    $$y^2=f(x),$$
where $f\in\FF_p[x]$ is of degree six and have no multiple roots \cite[see][chapter~1]{cassels}.

The set of principal divisors $\mathcal{P}(C)$ on $C$ constitutes a subgroup of the degree 0 divisors $\Div_0(C)$. The
Jacobian $\jac{C}$ of $C$ is defined as the quotient
    $$\jac{C}=\Div_0(C)/\mathcal{P}(C).$$
Let $\ell\neq p$ be a prime number. The $\ell^n$-torsion subgroup $\jac{C}[\ell^n]<\jac{C}$ of elements of order
dividing $\ell^n$ is then \cite[theorem~6, p.~109]{lang59}
    \begin{equation}\label{eq:J-struktur}
    \jac{C}[\ell^n]\simeq\ZZ/\ell^n\ZZ\times\ZZ/\ell^n\ZZ\times\ZZ/\ell^n\ZZ\times\ZZ/\ell^n\ZZ,
    \end{equation}
i.e. $\jac{C}[\ell^n]$ is a $\ZZ/\ell^n\ZZ$-module of rank four.

The order of $p$ modulo $\ell$ plays an important role in cryptography.

\begin{defn}[Embedding degree]
Consider a prime number $\ell$ dividing the order of $\jac{C}(\FF_p)$, where $\ell$ is different from $p$. The embedding degree of $\jac{C}(\FF_p)$ with respect to $\ell$ is the least number $k$, such that $p^k\equiv 1\mod{\ell}$.
\end{defn}

An endomorphism $\varphi:\jac{C}\to\jac{C}$ induces a $\ZZ_\ell$-linear map
    $$\varphi_\ell:T_\ell(\jac{C})\to T_\ell(\jac{C})$$
on the $\ell$-adic Tate-module $T_\ell(\jac{C})$ of $\jac{C}$ \cite[chapter~VII, \S1]{lang59}. The map $\varphi_\ell$ is given by $\varphi$ as described in the following diagram:
    $$
    \xymatrix@C=40pt@R=40pt{
    \dots \ar[r]^(0.4){[\ell]} & \jac{C}[\ell^{n+1}] \ar[r]^{[\ell]} \ar[d]^{\varphi} & \jac{C}[\ell^{n}] \ar[r]^{[\ell]} \ar[d]^{\varphi} & \dots \\
    \dots \ar[r]^(0.4){[\ell]} & \jac{C}[\ell^{n+1}] \ar[r]^{[\ell]} & \jac{C}[\ell^{n}] \ar[r]^{[\ell]} & \dots \\
    }
    $$
Here, the horizontal maps $[\ell]$ are the multiplication-by-$\ell$ map. Hence, $\varphi$ is represented by a
matrix $M\in\Mat_{4\times 4}(\ZZ/\ell\ZZ)$ on $\jac{C}[\ell]$. Let $P(X)\in\ZZ[X]$ be the characteristic polynomial of $\varphi$
\cite[see][pp.~109--110]{lang59}, and let $P_M(X)\in(\ZZ/\ell\ZZ)[X]$ be the characteristic polynomial of the restriction of $\varphi$ to $\jac{C}[\ell]$. Then \cite[theorem~3, p.~186]{lang59}
    \begin{equation}\label{eq:KarPolKongruens}
    P(X)\equiv P_M(X)\mod{\ell}.
    \end{equation}

Since $C$ is defined over $\FF_p$, the mapping $(x,y)\mapsto (x^p,y^p)$ is an isogeny on~$C$. This isogeny induces the $p$-power Frobenius endo\-morphism $\varphi$ on the Jacobian $\jac{C}$. The characteristic polynomial $P(X)$ of
$\varphi$ is of degree four \cite[theorem~2, p.~140]{tate}, and by the definition of $P(X)$
\cite[see][pp.~109--110]{lang59},
    $$|\jac{C}(\FF_p)|=P(1),$$
i.e. the number of $\FF_p$-rational elements of the Jacobian is determined by $P(X)$.

\section{CM fields}\label{sec:CMfields}

An elliptic curve $E$ with $\ZZ\neq\End(E)$ is said to have \emph{complex multiplication}. Let $K$ be an ima\-ginary,
quadratic number field with ring of integers $\heltal{K}$. $K$ is a \emph{CM field}, and if
\mbox{$\End(E)\simeq\heltal{K}$}, then $E$ is said to have \emph{CM by $\heltal{K}$}. More generally a CM field is
defined as follows.

\begin{defn}[CM field]
A number field $K$ is a CM field, if $K$ is a totally imaginary, quadratic extension of a totally real number field
$K_0$.
\end{defn}

In this paper only CM fields of degree $[K:\QQ]=4$ are considered. Such a field is called a \emph{quartic} CM field.

\begin{rem}\label{rem:quarticCM}
Consider a quartic CM field $K$. Let $K_0=K\cap\RR$ be the real subfield of $K$. Then $K_0$ is a real, quadratic number
field, $K_0=\QQ(\sqrt{D})$. By a basic result on quadratic number fields, the ring of integers of $K_0$ is given by
$\heltal{K_0}=\ZZ+\xi\ZZ$, where
    $$
    \xi=\begin{cases}
    \frac{1+\sqrt{D}}{2}, & \textrm{if $D\equiv 1\mod{4}$,} \\
    \sqrt{D}, & \textrm{if $D\equiv 2,3\mod{4}$}.
    \end{cases}
    $$
Since $K$ is a totally imaginary, quadratic extension of $K_0$, a number $\eta\in K$ exists, such that $K=K_0(\eta)$,
$\eta^2\in K_0$. The number $\eta$ is totally imaginary, and we may assume that $\eta=i\eta_0$, $\eta_0\in\RR$.
Furthermore we may assume that $-\eta^2\in\heltal{K_0}$; so $\eta=i\sqrt{a+b\xi}$, where $a,b\in\ZZ$.
\end{rem}

Let $C$ be a hyperelliptic curve of genus two. Then $C$ is said to have CM by~$\heltal{K}$, if
$\End(C)\simeq\heltal{K}$. The structure of $K$ determines whether $C$ is irreducible. More precisely, the following
theorem holds.

\begin{teo}\label{teo:reducibel}
Let $C$ be a hyperelliptic curve of genus two with $\End(C)\simeq\heltal{K}$, where $K$ is a quartic CM field. Then $C$ is
reducible if, and only if, $K/\QQ$ is Galois with Galois group $\gal(K/\QQ)\simeq\ZZ/2\ZZ\times\ZZ/2\ZZ$.
\end{teo}

\begin{proof}
\cite[proposition~26, p.~61]{shi}.
\end{proof}

Theorem~\ref{teo:reducibel} motivates the following definition.

\begin{defn}[Primitive, quartic CM field]\label{def:CMfieldPrimitive}
A quartic CM field $K$ is called primitive if either $K/\QQ$ is not Galois, or $K/\QQ$ is Galois with cyclic Galois
group.
\end{defn}

The CM method for constructing curves of genus~two with prescribed endomorphism ring is described in detail by
\cite{weng03} and \cite{gaudry}. In short, the CM method is based on the construction of the class polynomials of a
primitive, quartic CM field $K$ with real subfield $K_0$ of class number $h(K_0)=1$. The prime number $p$ has to be
chosen such that $p=x\overline x$ for a number $x\in\heltal{K}$. By \cite{weng03} we may assume that
$x\in\heltal{K_0}+\eta\heltal{K_0}$.

\section{Properties of $\jac{C}(\FF_p)$}\label{sec:properties}

Consider a primitive, quartic CM field $K$ with real subfield $K_0$ of class number $h(K_0)=1$, and let $p$ be an
uneven prime number such that $p=x\overline x$ for a number $x\in\heltal{K_0}+\eta\heltal{K_0}$. The main result of
this paper, given by the following theorem, concerns a curve of genus two with $\heltal{K}$ as endomorphism ring.

{\samepage
\begin{teo}\label{teo:ed=1}
With the notation as in remark~\ref{rem:quarticCM}, let $C$ be a hyperelliptic curve of genus two defined over $\FF_p$
with $\End(C)\simeq\heltal{K}$. Assume that the $p$-power Frobenius under this isomorphism is given by the number
$\omega=c_1+c_2\xi+(c_3+c_4\xi)\eta$, where $c_i\in\ZZ$. Consider a prime number $\ell\mid|\jac{C}(\FF_p)|$ with
$\ell\neq p$, $\ell\nmid D$ and $\ell\nmid c_2$. Assume that the $\ell$-Sylow subgroup of $\jac{C}(\FF_p)$ is not
cyclic. Then $p\equiv 1\mod{\ell}$, i.e. the embedding degree of $\jac{C}(\FF_p)$ with respect to $\ell$ is one.
\end{teo}
}

\begin{proof}
Consider a prime number $\ell\mid |\jac{C}(\FF_p)|$ with $\ell\nmid pc_2D$. If $\ell=2$, then obviously $p\equiv
1\mod{\ell}$. Hence we may assume that $\ell\neq 2$. Assume that the $\ell$-Sylow subgroup $S$ of $\jac{C}(\FF_p)$ is
not cyclic. Then $S$ contains a subgroup $U\simeq(\ZZ/\ell\ZZ)^2$. So
    $$(\ZZ/\ell\ZZ)^2<\jac{C}(\FF_p)[\ell]<\jac{C}[\ell].$$
Let $\{e_1,e_2\}\subseteq\jac{C}(\FF_p)$ be a basis of $(\ZZ/\ell\ZZ)^2$. Expand by the isomorphism~\eqref{eq:J-struktur} this set to a basis $\{e_1,e_2,f_1,f_2\}$ of $\jac{C}[\ell]$. It then follows that $1$ is
an eigenvalue of the Frobenius with eigenvectors $e_1$ and $e_2$, i.e. $1$ is an eigenvalue of multiplicity at least two.

First we assume that $D\equiv 2,3\mod{\ell}$. Let $P(X)$ be the characteristic polynomial of the Frobenius. Since the
conjugates of $\omega$ are given by $\omega_1=\omega$, $\omega_2=\overline{\omega}_1$, $\omega_3$ and
$\omega_4=\overline{\omega}_3$, where
    $$\omega_3 = c_1-c_2\sqrt{D}+i(c_3-c_4\sqrt{D})\sqrt{a-b\sqrt{D}},$$
it follows that
    $$
    P(X) = \prod_{i=1}^4(X-\omega_i)
         = X^4-4c_1X^3+(2p+4(c_1^2-c_2^2D))X^2-4c_1pX+p^2.
    $$
Since $1$ is an eigenvalue of the Frobenius of multiplicity at least two, the characteristic polynomial $P(X)$ is
divisible by $(X-1)^2$ modulo $\ell$. Now,
    $$P(X)=Q(X)\cdot(X-1)^2+R(X),$$
where
    \begin{align*}
    R(X) = {} & 4(1-3c_1-(c_1-1)p+2(c_1^2-c_2^2D))X \\
              & +p^2-2p-4(c_1^2-c_2^2D)+8c_1-3.
    \end{align*}
Since $R(X)\equiv 0\mod{\ell}$, it follows that
    \begin{equation}\label{eq:*}
    1-3c_1-(c_1-1)p+2(c_1^2-c_2^2D) \equiv 0 \mod{\ell}.
    \end{equation}
Since $|\jac{C}(\FF_p)|=P(1)$, we know that
    \begin{equation}\label{eq:***}
    (p+1)^2-4c_1(p+1)+4(c_1^2-c_2^2D)\equiv 0\mod{\ell}.
    \end{equation}
By equation~\eqref{eq:*} we see that $4(c_1^2-c_2^2D)\equiv
2(c_1-1)p-2+6c_1\mod{\ell}$. Substituting this into equation~\eqref{eq:***} we get
    $$(p+1)^2-4c_1(p+1)+2(c_1-1)p-2+6c_1\equiv 0\mod{\ell};$$
so either $p\equiv 1\mod{\ell}$ or $p\equiv 2c_1-1\mod{\ell}$. Assume $p\equiv 2c_1-1\mod{\ell}$. Then
    $$R(X)\equiv 4c_2^2D(-2X+1)\equiv 0\mod{\ell}.$$
Since $\ell\nmid 2c_2D$, this is a contradiction. So if $D\equiv 2,3\mod{4}$, then $p\equiv 1\mod{\ell}$.

Now consider the case $D\equiv 1\mod{4}$. We now have
    $$\omega_3 = c_1+c_2\frac{1-\sqrt{D}}{2}+i\left(c_3+c_4\frac{1-\sqrt{D}}{2}\right)\sqrt{a+b\frac{1-\sqrt{D}}{2}},$$
and it follows that the characteristic polynomial of the Frobenius is given by
    $$P(X)=X^4-2cX^3+(2p+c^2-c_2^2d)X^2-2pcX+p^2,$$
where $c=2c_1+c_2$. We see that $P(X)=Q(X)(X-1)^2+R(X)$, where
    \begin{align*}
    R(X) = {} & ((4-2c)p+2c^2-6c-2c_2^2D+4)X \\
              & +p^2-2p-3+4c-c^2+c_2^2D.
    \end{align*}
Since $R(X)\equiv 0\mod{\ell}$, it follows that
    \begin{equation}\label{eq:*1}
    p^2-2p-3+4c-c^2+c_2^2D \equiv 0 \mod{\ell},
    \end{equation}
and since $|\jac{C}(\FF_p)|=P(1)$, we know that
    \begin{equation}\label{eq:***1}
    (p+1)^2-2c(p+1)+c^2-c_2^2D \equiv 0 \mod{\ell}.
    \end{equation}
From equation~\eqref{eq:*1} and \eqref{eq:***1} it follows that
    $$p^2-cp+c-1 \equiv 0 \mod{\ell},$$
i.e. $p\equiv 1\mod{\ell}$ or $p\equiv c-1\mod{\ell}$. Assume $p\equiv c-1\mod{\ell}$. Then
    $$R(X)\equiv c_2^2D(-2X+1)\equiv 0\mod{\ell},$$
again a contradiction. So if $D\equiv 1\mod{4}$, then $p\equiv 1\mod{\ell}$.
\end{proof}

Consider the case $\ell\mid c_2$. Then the characteristic polynomial of the Frobenius modulo $\ell$ is given by
    $$P(X)\equiv (X^2-2c_1X+p)^2 \mod{\ell},$$
independently of the remainder of $D$ modulo $4$. Observe that
    $$X^2-2c_1X+p = (X+1-2c_1)(X-1)+p-2c_1+1.$$
Hence, $p\equiv 2c_1-1\mod{\ell}$, i.e.
    $$P(X)\equiv (X-1)^2(X-p)^2\mod{\ell}.$$
So the following theorem holds.

\begin{teo}\label{teo:c2}
With the notation as in remark~\ref{rem:quarticCM}, let $C$ be a hyperelliptic curve of genus two defined over $\FF_p$
with $\End(C)\simeq\heltal{K}$. Assume that the $p$-power Frobenius under this isomorphism is given by the number
$\omega=c_1+c_2\xi+(c_3+c_4\xi)\eta$, where $c_i\in\ZZ$. Consider a prime number $\ell\mid|\jac{C}(\FF_p)|$ with
$\ell\neq p$, $\ell\mid c_2$. Assume that the $\ell$-Sylow subgroup of $\jac{C}(\FF_p)$ is not cyclic. Then either
    \begin{enumerate}
    \item $\jac{C}(\FF_p)[\ell]\simeq(\ZZ/\ell\ZZ)^2$, or
    \item $p\equiv 1\mod{\ell}$ and $\jac{C}(\FF_p)[\ell]=\jac{C}[\ell]$.
    \end{enumerate}
\end{teo}

\begin{proof}
If $p\not\equiv 1\mod{\ell}$, then $1$ is not an eigenvalue of the Frobenius of multiplicity three, i.e.
$\jac{C}(\FF_p)[\ell]\simeq(\ZZ/\ell\ZZ)^2$. If $p\equiv 1\mod{\ell}$, then $1$ is an eigenvalue of the Frobenius of
multiplicity four, i.e. $\jac{C}(\FF_p)[\ell]=\jac{C}[\ell]$.
\end{proof}

\section{Applications}

Let $C$ be a hyperelliptic curve of genus two defined over $\FF_p$ with
$\End(C)\simeq\heltal{K}$. Write
    \begin{equation}\label{eq:jac(Fp)}
    \jac{C}(\FF_p)\simeq\ZZ/n_1\ZZ\times\ZZ/n_2\ZZ\times\ZZ/n_3\ZZ\times\ZZ/n_4\ZZ,
    \end{equation}
where $n_i\mid n_{i+1}$ and $n_2\mid p-1$ \cite[see][proposition~5.78, p.~111]{hhec}. We recall the following result on
the prime divisors of the number~$n_2$.

{\samepage
\begin{teo}\label{teo}
With the notion as above, let $\ell\mid n_2$ be an odd prime number. Then $\ell\leq Q$, where
    \begin{align*}
    Q &= \max\{a,D,a^2-b^2D\}, \\
    \intertext{if $D\equiv 2,3\mod{4}$, and}
    Q &= \max\{a,D,4a(a+b)-b^2(D-1),aD+2b(D-1)\},
    \end{align*}
if $D\equiv 1\mod{4}$. If $\ell>D$, then $c_1\equiv 1\mod{\ell}$ and $c_2\equiv 0\mod{\ell}$.
\end{teo}

\begin{proof}
\cite{me07a}.
\end{proof}
}

Let the Frobenius be given by the number $\omega=c_1+c_2\xi+(c_3+c_4\xi)\eta$, $c_i\in\ZZ$,
and consider a prime number $\ell\mid|\jac{C}(\FF_p)|$, $\ell\neq p$.

\begin{cor}
If $\ell\nmid c_2$ and $\ell>Q$, then the $\ell$-Sylow subgroup $S$ of $\jac{C}(\FF_p)$ is either of rank
two and $p\equiv 1\mod{\ell}$, or $S$ is cyclic.
\end{cor}

By \cite{me07b}, if $p\equiv 1\mod{\ell}$, then there exists an efficient, probabilistic algorithm to determine
generators of the $\ell$-Sylow subgroup of $\jac{C}(\FF_p)$. Hence the following corollary holds.

\begin{cor}
If $\ell\nmid D$ and $\ell\nmid c_2$, then there exists an efficient, probabilistic algorithm to
determine generators of the $\ell$-Sylow subgroup $S$ of $\jac{C}(\FF_p)$.
\end{cor}

\begin{proof}
If $p\equiv 1\mod{\ell}$, then the corollary is given by \cite{me07b}. If $p\not\equiv 1\mod{\ell}$, then $S$ is
cyclic by theorem~\ref{teo:ed=1}. Assume $|S|=\ell^n$. Then $S$ has $\ell^n-\ell^{n-1}$ elements of order $\ell^n$.
Hence the probability that a random element $\sigma\in S$ generates $S$ is $1-\ell^{-1}$, and choosing random elements
$\sigma\in S$ until an element of order $\ell^n$ is found will be an efficient, probabilistic algorithm to determine
generators of $S$.
\end{proof}

\section{Acknowledgement}

I would like to thank my supervisor Johan P. Hansen for inspiration on theorem~\ref{teo:c2}.


\begin{thebibliography}{99}

    \bibitem[Atkin and Morain(1993)]{atkin-morain} \textsc{A.O.L. Atkin and F. Morain}. Elliptic curves and primality
    proving. \emph{Math. Comp.}, vol.~61, pp.~29--68, 1993.
    \bibitem[Cassels and Flynn(1996)]{cassels} \textsc{J.W.S. Cassels and E.V. Flynn}. \emph{Prolegomena to a Middlebrow Arithmetic of Curves of Genus $2$}. London Mathematical Society Lecture Note Series. Cambridge University Press, 1996.
    \bibitem[Frey and Lange(2006)]{hhec} \textsc{G. Frey and T. Lange}. Varieties over Special Fields. In H.~Cohen and G.~Frey, editors, \emph{Handbook of Elliptic and Hyperelliptic Curve Cryptography}, pp.~87--113. Chapman \& Hall/CRC, 2006.
    \bibitem[Gaudry \emph{et al}(2005)Pierrick Gaudry]{gaudry} \textsc{P. Gaudry, T. Houtmann, D. Kohel, C. Ritzenthaler and A. Weng}.
    The $p$-adic CM-Method for Genus $2$. 2005. \url{http://arxiv.org}.%\url{http://arxiv.org/abs/math/0503148}.
    \bibitem[Koblitz(1989)]{koblitz89} \textsc{N. Koblitz}. Hyperelliptic cryptosystems. \emph{J. Cryptology}, vol.~1,
    pp.~139--150, 1989.
    \bibitem[Lang(1959) Serge Lang]{lang59} \textsc{S. Lang}. \emph{Abelian Varieties}. Interscience, 1959.
    \bibitem[Ravnshøj(2007a)]{me07a} \textsc{C.R. Ravnshøj}. \emph{Large Cyclic Subgroups of Jacobians of Hyperelliptic Curves}. 2007a. \url{http://arxiv.org}.
    \bibitem[Ravnshøj(2007b)]{me07b} \textsc{C.R. Ravnshøj}. \emph{Generators of Jacobians of Hyperelliptic Curves}. 2007b. \url{http://arxiv.org}.
    \bibitem[Shimura(1998) Goro Shimura]{shi} \textsc{G. Shimura}. \emph{Abelian Varieties with Complex Multiplication and Modular Functions}.
    Princeton University Press, 1998.
    \bibitem[Spallek(1994)]{spallek} \textsc{A.-M. Spallek}. \emph{Kurven vom Geschlecht $2$ und ihre Anwendung in
    Public-Key-Kryptosystemen}. Ph.D. thesis, Institut für Experimentelle Mathe\-matik, Universität GH Essen, 1994.
    \bibitem[Tate(1966)]{tate} \textsc{J. Tate}. Endomorphisms of abelian varieties over finite fields. \emph{Invent. Math.}, vol.~2, pp.~134--144, 1966.
    \bibitem[Weng(2003)]{weng03} \textsc{A. Weng}. Constructing hyperelliptic curves of genus~$2$ suitable for
    crypto\-graphy. \emph{Math. Comp.}, vol.~72, pp.~435--458, 2003.
\end{thebibliography}
\end{document}